\documentclass{elsart3p}%{elsart}

\usepackage{amsmath,amssymb,amsfonts} % Typical maths resource packages
\usepackage{pstricks,subfigure,graphics}
\usepackage{amsbsy}
\usepackage{verbatim}
\newcommand{\afrac}[2]{\frac{%
  \raisebox{.25ex}{$\scriptstyle #1$}}
  {\raisebox{.25ex}{$\scriptstyle #2$}}}
\newcommand{\pfrac}[2]{\frac{%
  \raisebox{.25ex}{$\scriptstyle #1$}}{#2}}
\def\ch{{Christoffel\ }}
\def\mkf{{Markoff\ }}

\def\N{{\mathbb N}}
\def\Z{{\mathbb Z}}

\renewcommand{\tilde}{\widetilde}
\newcommand{\X}{\{a,b\}}

\newcommand{\bs}{\boldsymbol{s}}
\newcommand{\bu}{\boldsymbol{u}}
\newcommand{\bv}{\boldsymbol{v}}

\usepackage{pstricks}
\psset{unit=\unitlength}
\psset{griddots=6,gridlabels=0,gridcolor=gray,subgriddiv=1}
\newcommand{\ChrSegment}[1][blue]{\psline[linewidth=0.02,linecolor=#1]}
\newcommand{\ChrPath}[1][red]{\psline[linewidth=0.05,linecolor=#1]}

\definecolor{faintred}{rgb}{1,.5,.5}

\journal{Information Processing Letters} 
\begin{document} 
\begin{frontmatter} 

\title{A note on the Markoff condition and central words}

\author{Amy Glen},
\ead{amy.glen@gmail.com}
\ead[url]{http://www.lacim.uqam.ca/$\scriptstyle\sim$glen}
\author{Aaron Lauve}, 
\corauth[cor]{Corresponding author.} 
%\ead{aaronlauve@gmail.com}
\ead[url]{http://www.lacim.uqam.ca/$\scriptstyle\sim$lauve}
\author{Franco V. Saliola\corauthref{cor}}
\ead{saliola@gmail.com}
\ead[url]{http://www.lacim.uqam.ca/$\scriptstyle\sim$saliola}

\address{%
	LaCIM, %\\ 
	Universit\'e du Qu\'ebec \`a Montr\'eal, %\\ 
	Case Postale 8888, succursale Centre-ville, %\\ 
	Montr\'eal (Qu\'ebec) H3C 3P8, %\\ 
	CANADA %
	} 
\thanks{Thanks are due to C. Reutenauer and J. Berstel, who introduced the authors to \ch words at the Centre des Recherches Math{\'e}matiques (Montr{\'e}al), March, 2007.}

\begin{abstract} 
We define {\em Markoff words} as certain factors appearing in bi-infinite words satisfying the {\em Markoff condition}. We prove that these words coincide with {\em central words}, yielding a new characterization of {\em \ch words}.
\end{abstract} 

\end{frontmatter}

%-------------------------------------------------------------------
% Begin DOCUMENT
%-------------------------------------------------------------------

\maketitle

%-------------------------------------------------------------------
% Begin SECTION
%-------------------------------------------------------------------
\section{Introduction}
In studying the minima of certain binary quad\-ratic forms $AX^2 + 2BXY + CY^2$, \mkf \cite{Mar:79,Mar:80} introduced a necessary condition that a bi-infinite word $\bs$ must satisfy in order that it represent the continued fraction expansions of the two roots of $AX^2 + 2BX + C$. Over an alphabet $\X$, his condition essentially states that each
factor $x\tilde m ab my$ occurring in $\bs$, where $\tilde m$ is the word $m$ read in reverse and $\{x,y\}=\{a,b\}$, has the property that $x=b$ and $y=a$. We call such words $m$ {\em \mkf words} in what follows. See Definition \ref{def:markoff}.

From  \cite{Reu:06} (see also \cite[pg.~30]{CusFla:89}), it is known that the bi-infinite words satisfying the \mkf condition are precisely the {\em balanced words} of Morse and Hedlund \cite{MorHed:40}. After the work of A.~de~Luca \cite{Luc:81,Luc:97}, we know that palindromes now play a `central'  role in the study of such words.
Here, we establish the following new characterization of a particular family of palindromes called {\em central words}.
\medskip

\begin{thm} \label{T:main2}
A word is a \mkf word if and only if it is a central word.
\end{thm}
\medskip

Central words hold a special place in the rich theory of \emph{Sturmian words} (e.g., see \cite[Chapter 2]{Lot:02}). For instance, it follows from the work of de~Luca and Mignosi \cite{Luc:97,LucMig:94} that central words coincide with the palindromic prefixes of standard Sturmian words.

%Central words hold a special place in the rich theory of {\em Sturmian words} (e.g., see \cite[Chapter 2]{Lot:02}). For instance, it follows from the work of  that the  

As an immediate consequence of Theorem \ref{T:main2}, we obtain a new characterization of {\em Christoffel words} in Corollary \ref{Cor:Christoffel}. Since the \mkf condition is relatively unknown, we discuss it and its relationship to \ch words at greater length in Section \ref{sec:ch-words}. 

%-------------------------------------------------------------------
% Begin SECTION
%-------------------------------------------------------------------
\section{The Markoff condition}
\label{sec:markoff}

Fix an alphabet $\X$. A finite sequence $a_1$, $a_2, \ldots, a_n$ of elements from $\X$ is called a \emph{word} of {\em length} $n$ and is written $w = a_1a_2\cdots a_n$. The length of $w$ is denoted by $|w|$ and we denote by $|w|_a$ (resp.~$|w|_b$) the number of occurrences of the letter $a$ (resp.~$b$) in $w$. %The {\em empty word} $\epsilon$ has length $0$. 

A {\em right-infinite} (resp.~{\em left-infinite}, {\em bi-infinite}) word over $\{a,b\}$ is a sequence indexed by $\N^+$ (resp.~$\Z\setminus\N^+$, $\Z$) with values in $\{a,b\}$. For instance, a left-infinite word is represented by $\bu = \cdots a_{-2}a_{-1}a_0$ and a right-infinite word by $\bv = a_1a_2a_3\cdots$, and their {\em concatenation} gives the bi-infinite word $\bu\bv = \cdots a_{-2}a_{-1}a_{0}a_1a_2a_3\cdots$. Infinite words are typically typed in boldface.

If $v = a_1a_2\cdots$ is a finite or a right-infinite word, then its \emph{reversal} $\tilde v$ is the word $\cdots a_2a_1$. Similarly, if $\bu$ is a left-infinite word, then its reversal is the right-infinite word $\tilde \bu$. We define the reversal of a bi-infinite word $\bs =
\cdots a_{-2} a_{-1} a_0 a_1 a_2 \cdots$ by $\tilde \bs = \cdots a_2 a_1 a_0 a_{-1} a_{-2} \cdots$. A finite word $w$ is a \emph{palindrome} if $w = \tilde w$. 

A {\em factor} of a finite or infinite word $w$ is a finite word $v$ such that $w = uvu'$ for some words $u$, $u'$. 
\medskip

\begin{defn}\label{def:markoff}~%
$(1)$ 
Suppose $\bs$ is a bi-infinite word on the alphabet $\X$. We say that $\bs$ satisfies the \textbf{\em Markoff condition} if for each factorization $\bs=\tilde \bu xy \bv$ with $\{x,y\}=\{a,b\}$, one has either $\bu = \bv$ or $\bu =my\bu'$ and $\bv=mx\bv'$ for some finite word $m$ (possibly empty) and right-infinite words $\bu'$,~$\bv'$. 
\newline$(2)$ 
A (finite) word $m$ is a \textbf{\em \mkf word} if there exists a bi-infinite word $\bs$ satisfying the \mkf condition with a factorization of the form $\bs = \tilde \bu y \tilde m xy m x \bv$, where $\{x,y\} = \X$. 
\end{defn}
\medskip

Note that a bi-infinite word $\bs$ satisfies the \mkf condition if and only if its reversal $\tilde \bs$ does, and  $\bs$ does not satisfy the \mkf condition if and only if $\bs$ or $\tilde \bs$ contains a factor of the form $a \tilde m ab m b$ for some finite word $m$.

Words $\bs$ satisfying the \mkf condition fall into four
classes: the periodic class; two aperiodic classes; and an ultimately
periodic class. See Section \ref{sec:ch-words}. An example of each
type appears below. 

\begin{gather}
\label{eq:M1} \cdots (aabab)(aabab)(aabab)(aabab) \cdots \\[0.5ex]
\label{eq:M2} % line is pi/4 x + e
%\cdots ababababaabababaababababaababababaab \cdots \\[0.5ex]
\cdots       abaabababaababababaabababa     \cdots \\[0.5ex]
\label{eq:M3} % line is pi/4 x
\cdots abaababababaabababaabababa \cdots \\[0.5ex]
\label{eq:M4} \cdots (baaa)(baaa)baab(aaab)(aaab) \cdots 
\end{gather}

%%%%%
% Type (M2)
%\begin{figure}[!ht]
%\begin{center}
%\setlength\unitlength{0.3cm}
%\psset{unit=\unitlength}
%\begin{picture}(20,17)(-10,-8)
%\psgrid(-10,-6)(10,11)
%\ChrSegment(-10,-5.282)(10,10.718)
%\ChrPath
%(-10,-6)(-9,-6)(-9,-5)(-8,-5)(-8,-4)%xyxy
%(-7,-4)(-7,-3)(-6,-3)(-6,-2)(-5,-2)%xyxyx
%(-4,-2)(-4,-1)(-3,-1)(-3,0)(-2,0)%xyxyx
%(-2,1)(-1,1)(0,1)(0,2)(1,2)%yxxyx
%(1,3)(2,3)(2,4)(3,4)(3,5)%yxyxy
%(4,5)(5,5)(5,6)(6,6)(6,7)%xxyxy
%(7,7)(7,8)(8,8)(8,9)(9,9)%xyxyx
%(10,9)(10,10)%xy
%\end{picture}
%\caption{The word following the line $\ell(x) = \frac{\pi}4x + e$.}
%\end{center}
%\end{figure}

%%%%%
% Type (M3)
%\begin{figure}[!ht]
%%\hrule
%\begin{center}
%\setlength\unitlength{0.3cm}
%\psset{unit=\unitlength}
%\begin{picture}(20,16)(-10,-8)
%\psgrid(-10,-8)(10,8)
%\ChrSegment(-10,-7.85397)(10,7.85397)
%\ChrPath(-10,-8)(-9,-8)(-8,-8)(-8,-7)(-7,-7)(-7,-6)(-6,-6)(-6,-5)(-5,-5)(-5,-4)(-4,-4)(-3,-4)(-3,-3)(-2,-3)(-2,-2)(-1,-2)(-1,-1)(0,-1)(0,0)(1,0)(2,0)(2,1)(3,1)(3,2)(4,2)(4,3)(5,3)(6,3)(6,4)(7,4)(7,5)(8,5)(8,6)(9,6)(9,7)(10,7)
%\end{picture}
%\caption{The word following the line $\ell(x) = \frac{\pi}4x$.}
%%\label{}
%\end{center}
%%\hrule
%\end{figure}

%See Section \ref{sec:ch-words} for details.%, including a prescription for how to continue the aperiodic examples \eqref{eq:M2} and \eqref{eq:M3}. 

%-------------------------------------------------------------------
% Begin SECTION
%-------------------------------------------------------------------
\section{The balanced property}
\label{sec:balanced}

Observe that the above examples of bi-infinite words are ``balanced'' in the following sense.  
\medskip

\begin{defn} \label{D:balance}
A finite or infinite word $w$ over $\{a,b\}$ is said to be \textbf{\em balanced} if for any two factors $u$, $v$ of $w$ with $|u| = |v|$, we have $\bigl||u|_{a} - |v|_{a}\bigr| \leq 1$ (or equivalently, $\bigl||u|_{b} - |v|_{b}\bigr| \leq 1$), i.e., the number of $a$'s (or $b$'s) in each of $u$ and $v$ differs by at most $1$.
\end{defn}
\medskip 

This notion dates back to the seminal work of Morse and Hedlund \cite{MorHed:40}. More recently, Reutenauer proved the equivalence between the Markoff condition and the above balanced property:
\medskip

\begin{prop} \label{P:balance} {\em \cite[Theorem 3.1]{Reu:06}} A bi-infinite word $\bs$ satisfies the Markoff condition if and only if $\bs$ is balanced. \qed
\end{prop} 
\medskip

In Section \ref{sec:ch-words} we recount Reutenauer's Theorem~6.1 in \cite{Reu:06}, which gives a refinement of the balanced property and of the \mkf condition yielding the four classes illustrated by \eqref{eq:M1}--\eqref{eq:M4}.

%-------------------------------------------------------------------
% Begin SECTION
%-------------------------------------------------------------------
\section{Central words} \label{S:central}

There exist several equivalent ways to define central words (see \cite[Chapter 2]{Lot:02}). Here we choose to use the following definition, as proved in \cite[Proposition~2.2.34]{Lot:02} using results from \cite{Luc:97,LucMig:94}.
%\Note{should we specify which result in \cite{LucMig:94}? No, several results in [5] lead to the given definition of central words.}
\medskip

\begin{defn} \label{D:central} A word $w$ over $\{a,b\}$ is \textbf{\em central} if and only if $awb$ and $bwa$ are balanced. 
\end{defn}
\medskip

The following fact is especially pertinent.
\medskip

\begin{lem} \label{L:central} Any central word is a palindrome.
\end{lem} 
\vspace{-10pt}

\begin{pf} If $w$ is central, then $awb$ and $bwa$ are balanced by Definition \ref{D:central}.  Arguing by contradiction, suppose $w$ is not a palindrome. There exist
words $u$, $v$, $z$ and letters $\{x,y\} = \{a,b\}$ such that $w = uxv = zy \tilde u$. But then 
\[
  xwy = xuxvy = xzy\tilde u y,
\]
and the factors $xux$ and $y\tilde u y$ contradict the balanced property of $xwy$. \qed 
\end{pf}
\vspace{-10pt}

\noindent{\it Note.} Lemma \ref{L:central} also appears under a different guise in \cite[Lemma~7]{LucMig:94}. Also see Corollary 2.2.9 in \cite{Lot:02}.
The above proof is easily adapted to show directly that \mkf words are palindromes.
\medskip

We are now ready to prove Theorem \ref{T:main2}: \emph{a word is a \mkf word if and only if it is a central word.} 
\vspace{-10pt}

\begin{pf*}{Proof of Theorem \ref{T:main2}}
Suppose $m$ is a Markoff word. Let $\bs$ be a bi-infinite word satisfying the Markoff condition for which $y \tilde m x y m x$ is a factor, where $\{x, y\} = \{a,b\}$. The reversal of this factor, namely $x\tilde m y x m y$, is a factor of $\tilde\bs$, which also satisfies the Markoff condition. Therefore, the words $amb$ and $bma$ are factors of bi-infinite words satisfying the Markoff condition, and hence are balanced by Proposition \ref{P:balance}. Thus $m$ is central, by Definition \ref{D:central}.

Conversely, suppose $m$ is a central word. Then $m$ is a palindrome by Lemma~\ref{L:central}, and moreover $amb= a\tilde m b$ is balanced (Definition~\ref{D:central}). Therefore the word $a\tilde m bamb$ is also balanced and it can be viewed as a factor of some bi-infinite word satisfying the \mkf condition by Proposition \ref{P:balance}---specifically, a bi-infinite word of the type represented in \eqref{eq:M1}, with $amb$ repeated bi-infinitely. Thus $m$ is a \mkf word by Definition~\ref{def:markoff}(2). \qed
%
%$amb = a\tilde m b$ and $bma = b\tilde m a$ are balanced words. Therefore the words $ambamb = a\tilde m b a m b$ and $bmabma = b\tilde m a b m a$ are also balanced, and each can be viewed as a factor of some bi-infinite word satisfying the \mkf condition by Proposition \ref{P:balance} (specifically, words of type \eqref{eq:M1}, with $amb$ or $bma$ repeated bi-infinitely). Thus $m$ is a \mkf word by Definition~\ref{def:markoff}(2). \qed
\end{pf*}
\vspace{-10pt}

An immediate corollary is a new characterization of {\ch words} (defined in the next section).\medskip

\begin{cor} \label{Cor:Christoffel}
A word $m$ is a \mkf word if and only if $amb$ is a \ch word.
\end{cor} 
\vspace{-10pt}

\begin{pf}
From \cite[Chapter 2]{Lot:02}, a finite word $amb$ is \ch word if and only if $m$ is a central word, i.e., a \mkf word  (by Theorem \ref{T:main2}). \qed
\end{pf}
%\vspace{-10pt}

%-------------------------------------------------------------------
% Begin SECTION
%-------------------------------------------------------------------
\section{\ch words}\label{sec:ch-words}

This section describes four classes of words satisfying the \mkf
condition and how they naturally coincide with four classes of
balanced words.

%Before defining \ch words, we return briefly to the \mkf condition. 
If a bi-infinite word $\bs$ satisfies the \mkf condition, then it falls
into exactly one of the following classes. 

\smallskip 
\noindent
Let $\{x,y\}=\{a,b\}$. \smallskip
\begin{enumerate}
\item[$(M_1)$] 
The lengths of the \mkf words $m$ occurring in $\bs$ are bounded and $\bs$ 
cannot be written as $\tilde\bu xy \bu$ for some word $\bu$. 
%There exists $N(\bs)\in\N$ such that every
%factorization $\bs = \tilde \bu xy \bv$ with $\{x,y\} = \{a,b\}$,
%implies $\bu = my\bu'$ and $\bv = mx v'$ for some finite word $m$
%with $|m|\leq N(\bs)$. 
\item[$(M_2)$] 
The lengths of the \mkf words $m$ occurring in $\bs$ are unbounded and
$\bs$ cannot be written as $\tilde\bu xy \bu$ for some word $\bu$.% with $\{x,y\}=\{a,b\}$.
\item[$(M_3)$] There is exactly one $j\in\Z$ such that $\bs$ has
the factorization $\bs=\tilde \bu s_js_{j+1} \bu$ with $s_j \neq s_{j+1}$.
\item[$(M_4)$] $\bs$ is not of type $(M_1)$--$(M_3)$.
\end{enumerate} \smallskip
(Equivalently, $\bs$ is in $(M_4)$ iff there exist at
least two $i\in\Z$ such that $\bs=\tilde
\bu s_is_{i+1} \bu$ with $s_i \neq s_{i+1}$.)

The four examples \eqref{eq:M1}--\eqref{eq:M4} in Section \ref{sec:markoff} correspond,
respectively, to the classes $(M_1)$--$(M_4)$ above.
We now turn to constructing words in each of the above classes.
To achieve this, we present a geometric construction of \ch
words, which allows for a description of balanced bi-infinite words.

\smallskip

Fix $p,q \in \N$, with $p$ and $q$ relatively prime.  Let $\mathcal P$ denote
the path in the integer lattice from $(0,0)$ to $(p,q)$ that satisfies: (i)
$\mathcal P$ lies below the line segment $\mathcal S$ which begins at the
origin and ends at $(p,q)$; and (ii) the region in the plane enclosed by
$\mathcal P$ and $\mathcal S$ contains no other points of $\Z \times \Z$
besides those of $\mathcal P$.

Each step in $\mathcal P$ moves from a point $(x,y) \in \Z \times \Z$
to either $(x+1,y)$ or $(x,y+1)$, so we get
a word $L(p,q)$ over the alphabet $\X$ by encoding steps of the
first type by the letter $a$ and steps of the second type by the
letter $b$. See Figure \ref{fig:ch word}. 

The word $L(p,q)$ is called the
\textbf{(lower) \ch word} of slope $\afrac qp$. The \emph{upper \ch
words} are defined analogously.
\begin{figure}[!ht]
\centering
\setlength\unitlength{0.53cm}
\psset{unit=\unitlength}
\subfigure%[$L(5,3)= aabaabab $]
{
%\scalebox{.95}{%
\begin{picture}(5,3)(2,1)
\psgrid(2,1)(7,4)
\ChrSegment(2,1)(7,4)
\ChrPath(2,1)(3,1)(4,1)(4,2)(5,2)(6,2)(6,3)(7,3)(7,4)
\small
\uput{0.13}[270](2.5,1){$a$}
\uput{0.13}[270](3.5,1){$a$}
\uput{0.13}[270](4.5,2){$a$}
\uput{0.13}[270](5.5,2){$a$}
\uput{0.13}[270](6.5,3){$a$}
\uput{0.13}[180](4,1.5){$b$}
\uput{0.13}[180](6,2.5){$b$}
\uput{0.13}[180](7,3.47){$b$}
\end{picture}
%}%
}
\hspace{1em}
\subfigure%[$U(5,3)=babaabaa$]
{
%\scalebox{.95}{%
\begin{picture}(5,3)
\psgrid(0,0)(5,3)
\ChrSegment(0,0)(5,3)
\ChrPath(0,0)(0,1)(1,1)(1,2)(2,2)(3,2)(3,3)(4,3)(5,3)
\small
\uput{0.13}[90](0.5,1){$a$}
\uput{0.13}[90](1.5,2){$a$}
\uput{0.13}[90](2.5,2){$a$}
\uput{0.13}[90](3.5,3){$a$}
\uput{0.13}[90](4.5,3){$a$}
\uput{0.13}[0](0,0.53){$b$}
\uput{0.13}[0](1,1.5){$b$}
\uput{0.13}[0](3,2.5){$b$}
\end{picture}
%}%
}
\caption{The lower and upper Christoffel words of slope $\frac35$
are $aabaabab$ and $babaabaa$, respectively.}
\label{fig:ch word}
\end{figure}

For an introduction to the beautiful theory of \ch words, see \cite[Chapter 2]{Lot:02} or \cite{BLRS:xx}.
\smallskip

If the line segment $\mathcal S$ (as defined above) is replaced by a line
$\ell$, then the construction produces balanced bi-infinite words. 
Moreover, all balanced bi-infinite words can be obtained by modifying this
construction; they fall naturally into the following four classes determined by
$\ell$.
\smallskip

%\Note{
%CR's paper has $B_2$ and $B_3$ reversed! 
%}
%
%%%%%
%M1 - \mkf words are bounded
%M2 - \mkf words are not bounded
%M3 - unique factorization s = \tilde{u} xy u
%M4 - several factorizations s = \tilde{u} xy u
%
%MH1 - periodic (rational slope)
%MH2 - irrational line with no lattice points!
%MH3 - irrational line through the origin!
%MH4 - ultimately periodic, but not periodic

\begin{enumerate}
\item[$(B_1)$]
$\ell(x) = \pfrac qpx$ is a line of rational slope $\pfrac qp$ (these
are the periodic balanced words, see \eqref{eq:M1}).
\item[$(B_2)$]
$\ell$ is a line of irrational slope that does not meet 
any point of $\Z \times \Z$ (in \eqref{eq:M2}, $\ell(x) = \frac\pi4 x + e$). 
\item[$(B_3)$]
$\ell(x) = \alpha x$ is a line of irrational slope 
meeting exactly one point of $\Z \times \Z$ (in \eqref{eq:M3}, $\ell(x)=\frac\pi4 x$). 
\item[$(B_4)$]
The balanced words not of type $(B_1)$--$(B_3)$.
\end{enumerate} \smallskip

Balanced bi-infinite words of type $(B_4)$, represented in \eqref{eq:M4}, 
are either of the form $\cdots xx y xx \cdots$
or $\cdots (ymx)(ymx)(ymy)(xmy)(xmy)\cdots$,
where $\{x,y\}=\{a,b\}$ and $m$ is a \mkf word.
%Balanced bi-infinite words of type $(B_4)$ are either of the form 
%$\cdots xxx y xxx \cdots$, where $\{x,y\}=\{a,b\}$,
%or of the form 
%$\cdots (xmy)(xmy)(xmx)(ymx)(ymx)\cdots$,
%where $m$ is a \mkf word.
Hence, it is possible to adapt the geometric construction above
to construct this class of balanced words also.
%$\ell$ is a piecewise linear curve of the form,
%\begin{gather*}
%\ell(x) = 
%\begin{cases}
%rx+1, & x \in (1,\infty), \\
%r(x-1), & x \in (-\infty,1],
%\end{cases}
%\end{gather*}
%where $r$ is a rational number. 
See Figure \ref{fig:M4}.

\begin{figure}[!ht]
%\hrule
\centering
\setlength{\unitlength}{0.53cm}
\psset{unit=\unitlength}
\begin{picture}(14,6)(-5,-2)
\psgrid(-5,-2)(9,4)
\ChrSegment(9,4)(1,1.3333)
\ChrSegment(-5,-2)(1,0)
\ChrPath(1,0)(0,0)(-1,0)(-2,0)(-2,-1)(-3,-1)(-4,-1)(-5,-1)(-5,-2)
\ChrPath(1,0)(1,1)(2,1)(3,1)(3,2)(4,2)(5,2)(6,2)(6,3)(7,3)(8,3)(9,3)(9,4)
\footnotesize
\uput{0.13}[0](-5,-1.5){$b$}
\uput{0.13}[90](-4.5,-1){$a$}
\uput{0.13}[90](-3.5,-1){$a$}
\uput{0.13}[90](-2.5,-1){$a$}
\uput{0.13}[0](-2,-0.5){$b$}
\uput{0.13}[90](-1.3,0){$a$}
\uput{0.13}[90](-0.5,0){$a$}
\uput{0.13}[90](0.5,0){$a$} 
\uput{0.13}[0](1,0.4){$b$} %the middle piece
\uput{0.13}[270](1.6,1){$a$} %the middle piece
\uput{0.13}[270](2.53,1){$a$}
\uput{0.13}[180](3,1.5){$b$}
\uput{0.13}[270](3.5,2){$a$}
\uput{0.13}[270](4.5,2){$a$}
\uput{0.13}[270](5.5,2){$a$}
\uput{0.13}[180](6,2.5){$b$}
\uput{0.13}[270](6.5,3){$a$}
\uput{0.13}[270](7.5,3){$a$}
\uput{0.13}[270](8.5,3){$a$}
\uput{0.13}[180](9,3.5){$b$}
\end{picture}
\caption{Constructing example \eqref{eq:M4} $\cdots(baaa)baab(aaab)\cdots$.}
\label{fig:M4}
%\hrule
\end{figure}
As shown in \cite{Reu:06}, classes $(B_1)$--$(B_4)$ are derived from Morse and Hedlund's description of balanced bi-infinite words \cite{MorHed:40} (also see Heinis \cite{Hei:01}).

\medskip

\begin{prop} {\em \cite[Theorem 6.1]{Reu:06}} For $1\leq i \leq 4$, one has the coincidences $(M_i) = (B_i)$.
\end{prop}
\medskip

In closing, we mention that \mkf was interested in words over the
alphabet $\{1,2\}$ that satisfy the \mkf condition. For these words, he
studied the continued fraction quantities 
\[
\lambda_i(\bs) = s_i + [0,s_{i+1},\cdots] + [0,s_{i-1},s_{i-2},\cdots] 
\]
and $\Lambda(\bs)=\sup_i \lambda_i(\bs)$. 
Reutenauer \cite[Theorem~7.2]{Reu:06} showed that classes
$(M_1)$--$(M_4)$ correspond, respectively, to those $\bs$ satisfying
the \mkf condition with: $\Lambda(\bs)<3$; $\lambda_i(\bs) < 3$ for all $i$
but $\Lambda(\bs)=3$; $\Lambda(\bs)=3=\lambda_i(\bs)$ for a unique $i\in \Z$;
$\Lambda(\bs)=3=\lambda_i(\bs)$ for at least two $i\in \Z$. 

The set $\bigl\{ \Lambda(\bs) \mid \textrm{$\bs$ is a bi-infinite word over $\N^+$}\bigr\}$, with none of the conditions on $\bs$ originally imposed by Markoff, has become known as the \textbf{Markoff spectrum}. Results and open questions concerning the \mkf spectrum may be found in \cite{CusFla:89}.
%\vspace{-5pt}

\end{document}